\documentclass[reqno]{article}
\usepackage{amsmath, amsfonts, amscd, eucal}
\usepackage{latexsym}
\usepackage{amssymb}
\newcommand{\be}{\begin{equation}}
      \newcommand{\ee}{\end{equation}}
      \newcommand{\ba}{\begin{eqnarray}}
       \newcommand{\ea}{\end{eqnarray}}
\newcommand{\ban}{\begin{eqnarray*}}
\newcommand{\ean}{\end{eqnarray*}}

\newcommand{\pt}{\partial}

 \renewcommand{\o}[2]{\frac{#1}{#2}}
\newcommand{\hf}{\o{1}{2}}

 \newcommand{\qed}{\hspace*{\fill}\rule{3mm}{3mm}\quad \vspace{.2cm}}
 \newcommand{\Pf}{\noindent {\bf Proof:} }

\newcommand{\tr}{{\rm tr} }

\newtheorem{theo}{Theorem}[section]

\begin{document}
\newtheorem{defn}[theo]{Definition}
\newtheorem{ques}[theo]{Question}
\newtheorem{lem}[theo]{Lemma}
\newtheorem{prop}[theo]{Proposition}
\newtheorem{coro}[theo]{Corollary}
\newtheorem{ex}[theo]{Example}
\newtheorem{note}[theo]{Note}
\newtheorem{conj}[theo]{Conjecture}
\newtheorem{remark}[theo]{Remark}

\title{The Intersection R-Torsion for Finite Cone}
\author{Xianzhe Dai\thanks {Math Dept, UCSB, Santa Barbara, CA 93106 \tt{Email:
dai@math.ucsb.edu}. Partially supported by NSF} \and
Xiaoling Huang\thanks{Math Dept, UCSB, Santa Barbara, CA 93106, currently at EMC
 Partially supported by NSF  } }
\date{}
\maketitle

\section{Introduction}

Torsion invariants were originally introduced in the 3-dimensional
setting by K. Reidemeister \cite{R} in 1935 who used them to give a
homeomorphism classification of 3-dimensional lens spaces. The
Reidemeister torsions (R-torsions for short) are defined using
linear algebra and combinatorial topology. The salient feature of
R-torsions is that it is not a homotopy invariant but rather a
simple homotopy invariant; hence a homeomorphism invariant as well.
From the index theoretic point of view, R-torsion is a secondary
invariant with respect to the Euler characteristic. For geometric
operators such as the Gauss-Bonnet and Dolbeault operator, the index
is the Euler characteristic of certain cohomology groups. If these
groups vanish, the Index Theorem has nothing to say, and secondary
geometric and topological invariant, e.g., R-torsion, appears.
R-torsions were generalized to arbitrary dimensions by W. Franz
\cite{F} and later studied by many authors (Cf. \cite{M}).

Analytic torsion (or Ray-Singer torsion), which is a certain
combinations of determinants of Hodge Laplacians on k-forms, is an
invariant of Riemannian manifolds defined by Ray and Singer
\cite{RS} as an analytic analog of R-torsions. The Ray-Singer conjecture,  that the analytic and Reidemeister torsion
agree on closed manifolds,  was proved independently by
 Cheeger \cite{c1} and M\"uller \cite{mu1} using
different techniques. Later, Vishik \cite{V} gave a cutting and
pasting proof from the viewpoint of topological quantum field theory, 
and Bismut and Zhang \cite{bz} used Witten deformation
to generalize it to arbitrary flat bundles, see also \cite{mu2} for the unimodular case.

Further significant work includes that of M\"uller \cite{mu2}, which
extended the theorem to unimodular representations, that of Bismut
and Zhang \cite{bz}, which treated general representations (in which
interesting secondary invariants come in), and that of
Burghelea-Friedlander-Kappeler-McDonald \cite{BFKM}, which dealt
with infinite dimensional representations.

It is a natural question wether the Ray-Singer
conjecture/Cheeger-M\"uller theorem extends to singular manifolds.
For manifolds with isolated conical singularity, both the R-torsion
and analytic torsion have been defined by Dar \cite{dar}, using
respectively, the intersection homology of Goresky-MacPherson
\cite{gm1,gm2} and Cheeger's theory of heat kernels for conical
singularity \cite{c2}. There are several possible approaches to this
question, among which the most natural one is to reduce the problem
to three parts. One concerns manifolds with boundary, for which the
question has been extensively studied \cite{c1,LR,lu,DF,BM}. The
second part would be a finite cone. The last part deals with the
Mayer-Vietories principle.

In this paper we concentrate on the intersection R-torsion side of
the story. (We note that there has been extensive work on the analytic torsion of
cones and conical manifolds, see, for example, \cite{s,ve,hs,mv}.) We will first study the intersection R-torsion of a
finite cone. Our main result expresses it as a combination of
determinants of the combinatorial Laplacian on the cross section of
the cone. We then study an analytic invariant which is obtained by
replacing the combinatorial Laplacian by the Hodge Laplacian.

More specifically, consider the finite cone $X=C(Y)$ with the cross
section $Y$ a closed $(n-1)$-dimensional manifold. Let
$I\tau^{\bar{p}}(X)$ denote the intersection R-torsion of $X$, where
$\bar{p}=(p_2, p_3, \cdots, p_n)$ is a given perversity. Then,

\begin{theo} Let $\Delta^{(c)}$ denote the combinatorial Laplacian of the cross section $Y$.
Then $$\ln I\tau^{\bar{p}}(X)= \sum_{p=0}^{n-p_n-1} (-1)^{p+1}p\,
\ln \det \Delta_p^{(c)} + (n -p_n) \sum_{p-p_n}^{n-1}
(-1)^{p+1}\ln \det \Delta_p^{(c)}. $$
\end{theo}

This leads us to an analogous analytic invariant for an even dimensional manifolds. Thus, let $Y$ be an even dimensional closed manifold with
$m=\dim Y$. Let $p$ be an integer such that $0\leq p \leq m-1$ ($p$
corresponds to $p_n$ which is determined by a given perversity).
Given an orthogonal representation $\rho: \, \pi_1(Y)
\longrightarrow O(N)$, one has an associated flat vector bundle
$E_{\rho}$ with compatible metric on $Y$. Let $\Delta_k$ be the
Laplacian acting on differential $k$ forms on $Y$ with coefficients
in $E_{\rho}$. Then we define

\[  \ln T_{p}(Y, \rho) = \frac{1}{2} \left[
\sum_{k=0}^{m-p} (-1)^{k+1} k\, \ln \det(\Delta_k) +
(m-p)\sum_{k=m-p+1}^m (-1)^{k+1} \ln \det(\Delta_k) \right] . \]

For $p=0$ this gives the usual analytic torsion which is trivial for dimensional reasons.
Other values of $p$ give nontrivial and more interesting analytic invariants that bear close resemblance to the so called Cheeger's half torsion, see \S 4.
To investigate what kind of invariant $\ln T_{p}(Y, \rho)$ defines,
we now look at its variation under metric change. Let $g(u)$ be a
family of Riemannian metrics on $Y$ and $\Delta_k(u)$ the
corresponding Laplacian (when there is no ambiguity we will often
write $\Delta_k$ instead of $\Delta_k(u)$.  Let $\dot{\star}=
d\star/du$ and $\alpha=\star^{-1} \dot{\star}$.  Denote by
$E_k(t)=e^{-t\Delta_k(u)}$ the heat kernel and let $E_k=E_k^{ex} +
E_k^{ce} + E_k^h$ denote the Hodge decomposition of $E_k$ into its
exact, coexact and harmonic parts. We have the following result
regarding the variation of $\ln T_{p}(Y, \rho)$.

\begin{theo}  The variation of $\ln T_p(Y, \rho)$ is given by
\ban \frac{d}{du} \ln T_p(Y, \rho) & = & \hf \sum_{k=0}^{m-p-1}
(-1)^{k+1} {\rm Tr}(P_{H^k} \alpha) + \hf \sum_{k=0}^{m-p-1}
(-1)^{k+1} {\rm LIM}_{t\rightarrow 0} {\rm
Tr}(e^{-t\Delta_k} \alpha)\\
& + & (-1)^{m-p+1} \hf {\rm LIM}_{t\rightarrow 0} {\rm
Tr}(E_{m-p}^{ex}(t) \alpha), \ean where $P_{H^k}$ denote the
projection onto the cohomology $H^k$ and ${\rm LIM}_{t\rightarrow 0}
{\rm Tr}(E_{m-p}^{ex}(t) \alpha)$ denotes the constant term in the
asymptotic expansion of ${\rm Tr}(E_{m-p}^{ex}(t) \alpha)$. 

In particular, for $p=m/2$, the variation of this invariant under conformal changes is local:
\[ \frac{d}{du} \ln T_{\frac{m}{2}}(Y, \rho)  =  \hf \sum_{k=0}^{\frac{m}{2}-1}
(-1)^{k+1} {\rm Tr}(P_{H^k} \alpha) + \hf \sum_{k=0}^{\frac{m}{2}-1}
(-1)^{k+1} {\rm LIM}_{t\rightarrow 0} {\rm
Tr}(e^{-t\Delta_k} \alpha). \]
\end{theo}

Finally, we examine the R-torsion of the Mayer-Vietoris sequence.
\begin{theo}
Assume that the Witt condition $H^{\frac{m}{2}}(Y)=0$ holds. Then
the R-torsion of the Mayer-Vietoris sequence in intersection
cohomology
$$ \cdots \longrightarrow IH^{q}_{(2)}(Y)\longrightarrow IH_{(2)}^{q+1}(X)\longrightarrow IH_{(2)}^{q+1}(M)\oplus
IH_{(2)}^{q+1}(C(Y))\longrightarrow IH_{(2)}^{q+1}(Y)\longrightarrow
\cdots
$$
is equal to the R-torsion of the truncated exact sequence of the
pair $(M,Y)$
$$ 0 \longrightarrow H^{\frac{m}{2}+1}(M,Y) \longrightarrow H^{\frac{m}{2}+1}(M)\longrightarrow
H^{\frac{m}{2}+1}(Y)\longrightarrow H^{\frac{m}{2}+2}(M,Y)
\longrightarrow \cdots $$
\end{theo}

The rest of the paper is organized as follows. In the next section, we review the definition of R-torsion and recall the intersection cohomology of Goresky-MacPherson, leading to the definition of intersection R-torsion. In \S 3, we compute the intersection R-torsion of a finite cone. The analytic analog is studied in \S4. Finally, \S 5 deals with the R-torsion of the Mayer-Vietoris sequence.

The first author would like to acknowledge very useful communications with Jeff Cheeger, Rafe Mazzeo and Boris Vertman.

\section{The definition of Intersection R-torsion}
 We briefly recall the definition and characteristic
properties of R-torsion. Roughly
speaking, the R-torsion measures to what extent the boundary map of
a chain complex can be made to preserve a preferred volume element.
Let $C$ be a real vector space of dimension $n$ and let $b=(b_1,
\cdots , b_n), \ c=(c_1,\cdots , c_n)$ be two different bases for
$C$, Then $c_i=a_{ij}b_j$ and $(a_{ij}) \in GL(n,R)$. We denote the volume change between two bases
$det(a_{ij} )$ by $[c/b]$.

Let $(C, \partial):  0 \rightarrow C_n
\stackrel{\partial_n}{\rightarrow} C_{n-1}
\stackrel{\partial_{n-1}}{\rightarrow} \cdots C_1
\stackrel{\partial_1}{\rightarrow} C_0 \rightarrow 0$ be a chain
complex of real vector spaces. Let $c_i$ be a preferred basis for
$C$ and $h$ a preferred basis for the homology group $H^*(C)$.
Denote by $B_i$ the image of the boundary map $\partial_{i+1}: \
C_{i+1} \rightarrow C_i$ and $Z_i$ its kernel. We choose a basis
$b_i$ for $B_i$, which lifts to linearly independent set
$\tilde{b}_{i} \in C_{i+1}$, i.e. $\partial \tilde{b}_{i}=b_i$.
Using the inclusions $0 \subset B_i \subset Z_i \subset C_i$ where
$Z_i/B_i \equiv H_i,\ C_i/Z_i \equiv B_{i-1}$ we see that $ b_i,
h_i, \tilde{b}_{i-1}$ combine to give a new basis for $C_i$. The
R-torsion of the chain complex is the real number $\tau(c, h)$
defined by \be \label{rtocc} \ln \tau(c, h)= \sum_{i=0}^n (-1)^i \ln
|[b_i h_i \tilde{b}_{i-1}/c_i]|. \ee

The R-torsion $\tau(c, h)$ does not depend on the choice of $b_i,
\tilde{b}_{i-1}$, but it depends on the preferred bases $c_i, h_i$.
In fact, it depends only on the volume elements determined by these
preferred bases. More precisely,

\be \label{bcfrt1} \ln \tau(c', h)=\ln \tau(c, h) + \sum_{i=0}^n (-1)^i \ln
|[c_i/c'_i]|, \ee
and 

\be \label{bcfrt2} \ln \tau(c, h')=\ln \tau(c, h) + \sum_{i=0}^n (-1)^i \ln
|[h'_i/h_i]|. \ee

When the preferred basis of the homology is chosen according to the
preferred basis of the chain complex, there is an elegant
representation of the R-torsion in terms of the combinatorial
Laplacian. The choice of a preferred basis for each $C_i$ represents
$\partial_{i}: \ C_{i} \rightarrow C_{i-1}$ as a real matrix. Let
$\partial_{i}^*: \ C_{i-1} \rightarrow C_i$ be the transpose matrix.
The combinatorial Laplacian is $\Delta_i^{(c)}=\partial_{i+1}
\partial_{i+1}^* + \partial_i^* \partial_i:\ C_i \rightarrow C_i$.
By the finite dimensional Hodge theory, $\ker \Delta_i^{(c)} \cong
H_i(C, \partial)$. Effectively, the preferred basis $c_i$ determines
an inner product on $C_i$ in which $c_i$ becomes orthonormal. If we
choose the preferred basis $h$ on $H_i(C,
\partial)$ to correspond to an orthonormal basis of $\ker
\Delta_i^{(c)} \subset C_i$, then, \be \label{rtitcl} \ln \tau(c,
h)= \frac{1}{2} \sum_{i=0}^n (-1)^{i+1} i \log \det \Delta_i^{(c)}.
\ee

Now if $K$ be a finite CW complex, consider $\tilde{K}$ the
universal covering complex of $K$. The fundamental group $\pi$ of
$K$ acts on $\tilde{K}$ as the group of covering transformations.
This action makes $C(\tilde{K})$, the cellular chain complex
associated with $\tilde{K}$,  a free $\mathbb R \pi$-module
generated by the cells $e_i$ of the complex $K$. We pick a preferred
basis for $C_i(\tilde{K} )$ coming from the $i$-cells of $K$,
denoted $(e_i^1, e_i^2, \cdots , e_i^{k_i})$.

Let $\epsilon: \ \pi \longrightarrow O(n)$ be an orthogonal
representation of the fundamental group. Then one can construct a
chain complex of real vector spaces by setting $C_i(K,\epsilon)=
C_i(\tilde{K})\otimes_{\mathbb R \pi} R^n$. We have a preferred
choice of basis for each vector space $C_i(K, \epsilon)$ given by
$e_i^j \otimes x_k$ where $x_k$ is an orthonormal basis for $R^n$.
With a choice of preferred basis $h$ in homology, the torsion of the
complex of real vector spaces $C_i(K, \epsilon)$ is a real number
and will be denoted $\tau(K, \epsilon, h)$.

 The R-torsion is a
combinatorial invariant i.e. invariant under subdivision of $K$. It
is a topological invariant when the chain complex is acyclic.

The R-torsion of a closed manifold $M$ is the R-torsion of the cell
complex determined by a cell structure of $M$. In this case, the
preferred base for the homology is obtained via Hodge theory through
an orthonormal basis of the harmonic forms. With this choice of
preferred basis in homology it was shown in celebrated work of
Cheeger \cite{c1} and M\"uller \cite{mu1} that $\tau(M, \epsilon)$
equals the so called analytic torsion (Ray-Singer conjecture).

The intersection R-torsion is defined for pseudomanifolds by Dar
\cite{dar} using the intersection homology theory of
Goresky-MacPherson. We recall the basic facts of Intersection
Homology Theory.

A pseudomanifold $X$ of dimension $n$ is a compact PL space for
which there exists a closed subspace $Z$ with dimension $Z$ $\leq
n-2$ such that $X-Z$ is an n-dimensional oriented manifold which
is dense in $X$. A stratification of a pseudomanifold is a
filtration by closed subspaces \[ X=X_n= X_{n-1} \supset X_{n-2}=Z
\supset
 \cdots \supset X_1 \supset X_0\] such that for each point $p\in X_i-X_{i-1}$
there is a filtered space $V = V_n \supset V_{n-1}\supset \cdots V_i
=$a point and a mapping $V\times B^i \rightarrow X$ which takes
$V_j\times B^i$ (PL) homeomorphically to a neighborhood of p in
$X_j$.  $X_i-X_{i-1}$ is an $i$-dimensional manifold called the
$i$-dimensional stratum. Every pseudomanifold admits a
stratification.

We will actually be working with pseudomanifold with boundaries.
However, in our situation, the boundaries do not intersect with
the singular strata. Hence the above discussions can be easily
modified to adapt to the corresponding situation.

The space of geometric chains $C_*(X)$ is the collection of all
simplicial chains with respect to some triangulation where one
identifies the two chains if their images coincide under some common
subdivision. The intersection homology theory is obtained by
restricting to only {\em allowable} chains, described by the so
called perversity.

A perversity is a sequence of integers $\bar{p}=(p_2, p_3,
\cdots,p_n)$ such that $p_2=0$ and $p_{k+1} =p_k \ {\rm or} \ p_{k+
1}$. If $i$ is an integer and $\bar{p}$ is a perversity, a subspace
$Y\subset X$ is $(\bar{p},i)$ allowable if $\dim(Y)\leq i$ and
$\dim(Y\cap X_{n-k}) \leq i-k +p_k$ for $k\geq 2$. In other words,
$p_k$ describes how much $X$ is allowed to deviate from intersecting
the stratum $X_{n-k}$ transversally. The intersection chains
$IC_i^{\bar{p}}(X)$ is the subspace of $C_i(X)$ consisting of those
chains $\xi$ such that $|\xi|$ is $(\bar{p}, i)$ allowable and
$|\partial \xi|$ is $(\bar{p}, i-1)$ allowable. The $i$-th
Intersection Homology Group of perversity $\bar{p}$,
$IH_i^{\bar{p}}(X)$ is the $i$-th homology group of the chain
complex $IC_*^{\bar{p}}(X)$.

The intersection chain complex as we defined is not finitely
generated. In order to define the Intersection R-torsion we need to
work with finitely generated chain groups. To do this one uses the
basic sets $R_i^{\bar{p}}$.

 Let $X$ be a pseudomanifold with a fixed
stratification. Let $T$ be a triangulation of $X$ subordinate to the
stratification i.e. such that each $X_k$ is a subcomplex of $T$.
Define $R_i^{\bar{p}}$ be the subcomplex of $T'$, the first
barycentric subdivision of $T$, consisting of all simplices which
are $(\bar{p}, i)$ allowable.

Let ${\cal R}^{\bar{p}}(X)$ be the chain complex whose $i$-th chain
group consists of simplicial chains $e_i$ such that $|e_i| \in
R_i^{\bar{p}}$ and $|\partial e_i| \in R_{i-1}^{\bar{p}}$.  It is a
free abelian group generated by finitely many chains $\{e_i^j \}$.
The homology group $H_i({\cal R}^{\bar{p}}(X))$ is canonically
isomorphic to $IH_i^{\bar{p}}(X)$.

Let $\tilde{X}$ be the universal covering complex of $X$. Then the
chain complex ${\cal R}^{\bar{p}}(\tilde{X})$ is a free $\mathbb R
\pi$-module generated by the lifts of the chains $\{e_i^j \}$. If
$\epsilon:\ \pi \rightarrow O(n)$ be an orthogonal representation
one obtain a chain complex of real vector spaces ${\cal
R}^{\bar{p}}(X, \epsilon)= {\cal R}^{\bar{p}}(\tilde{X})
\otimes_{\mathbb R \pi} \mathbb R^n$ with a preferred basis given by
$\{ e_i^j \otimes x_k \}$ where $x_k$ is an orthonormal basis for
$\mathbb R^n$.

The intersection R-torsion of $X$ is then defined to be the torsion
of the chain complex ${\cal R}^{\bar{p}}(X, \epsilon)$, provided a
preferred basis in homology is chosen. Dar \cite{dar} proved that
the intersection R-torsion is a combinatorial invariant and
independent of the stratification.

\section{Intersection R-torsion of finite cone}

In this section we restrict ourself to the finite cone.  Let
$X=C(Y)$ be a  finite cone with $dim X=n$, where the cross section
$Y$ is a closed manifold.  We will also write $X=w\ast Y$ with $w$
the cone tip.
\newline

If $\sigma=[a_{0}, \cdots, a_{p}]$ is an oriented simplex of $Y$,
then $[w, \sigma]=[w, a_{0}, \cdots, a_{p} ]=w\ast \sigma$ is an
oriented simplex of $w\ast Y$. Similarly, if $\eta=\sum n_{i}\sigma
_{i}$ is a $p$-chain of $Y$, then $[w, \eta]=\sum
n_{i}[w,\sigma_{i}]$, and
\begin{align*}
\partial[w, \eta]=\begin{cases}
\eta-w & \mathrm{dim} \, \eta =0\\
\eta-[w , \partial {\eta}] &\mathrm{ dim}\, \eta >0
\end{cases}
\end{align*}
We have the following results:
\begin{lem} If $Z_{p}$ is a $p$-cycle of $X$ for $p\geq 1$, then $Z_{p}=\partial [w,
C_p]$ for some $p$-chain $C_p$ of $Y$.
\end{lem}
\Pf Write $Z_{p}=C_{p}+[w, D_{p-1}]$, where $C_{p}$ and $D_{p-1}$
are both carried by $Y$. Then by the above observation and using
that $Z_p$ is closed, we have $Z_{p}=\partial [w, C_p]$. \qed
\newline

Of course, this is compatible with the well known fact that
\begin{align*}
H_{p}(X)=\begin{cases}
0 &p\geq1\\
\mathrm{Z} &p=0
\end{cases}
\end{align*}

Now let $\bar{p}$ be a perversity. Since $X$ has only strata of
dimension $n$ and $0$, the intersection chains and homology will
only depend on $p_n$. In fact, we have
\begin{lem} \label{icc} The intersection chains of $X$ are given by
\[ IC_i^{\bar{p}}(X)=\left\{ \begin{array}{ll} C_i(Y), & i<n-p_n, \\
\{\ \xi \in C_i(X) \ | \ \partial \xi \in C_{i-1}(Y)\ \}, &
i=n-p_n,
\\
C_i(X), & i
> n-p_n.
\end{array} \right.
\]
In particular,
\be \label{icoc} IH_i^{\bar{p}}(X)=\left\{ \begin{array}{ll} H_i(Y), & i<n-p_n-1, \\
0, & i\geq n-p_n-1.
\end{array} \right.
\ee
\end{lem}

From now on we will sometimes suppress the superscript $\bar{p}$ here when there is no confusion.  Let $h_{p}(Y)=\{ h_{1}^{p}(Y), \cdots, h_{j_{p}}^{p}(Y)\}$ be a preferred basis for
$H_{p}(Y)$. Any triangulation of $Y$ gives rise to a preferred basis for $C_p(Y)$.  At first we will assume that the basis of $IH_i^{\bar{p}}(X)$ is chosen to be compatible with those of  $H_{p}(Y)$ with respect to  (\ref{icoc}).

\begin{theo} \label{mt} Let $\Delta^{(c)}$ denote the combinatorial Laplacian of the cross section $Y$.
Then $$\ln I\tau^{\bar{p}}(X)= \sum_{p=0}^{n-p_n-1} (-1)^{p+1}p\,
\ln \det \Delta_p^{(c)} + (n -p_n) \sum_{p-p_n}^{n-1}
(-1)^{p+1}\ln \det \Delta_p^{(c)}. $$
\end{theo}

\Pf The intersection R-torsion is defined in terms of the chain
complex
\begin{align}\label{eqn:4.1} \cdots \longrightarrow IC_{p+1}(X)\longrightarrow IC_{p}(X)\longrightarrow IC_{p-1}(X)
\longrightarrow \cdots.
\end{align}
We examine the terms of this complex according to their degrees.
\newline

\noindent Case 0. $p=n$\\
\newline
In this case, $IC_{n}(X)=C_{n}(X)$, so we only need to consider the usual chains of $X$.\\
\newline
Let $c_{n-1}(Y)=\{ \sigma_{1}^{n-1}(Y), \cdots,
\sigma_{i_{n-1}}^{n-1}(Y)\}$ be the preferred basis of
$(n-1)$-chains of $Y$.  Then \\
\{$[w, c_{n-1}(Y)]$\} is the preferred basis of $C_{n}(X)$. Choose a
basis $b_{p}(Y)=\{b_{1}^{p}(Y), \cdots, b_{k_{p}}^{p}(Y)\}$ for
$B_{p}(Y)$, and their lifts $\tilde{b}_{p}(Y)=\{
\tilde{b}_{1}^{p}(Y), \cdots, \tilde{b}_{k_{p}}^{p}(Y)\}$, and
$h_{p}(Y)=\{ h_{1}^{p}(Y), \cdots, h_{j_{p}}^{p}(Y)\}$ the basis for
$H_{p}(Y)$. Then by the fact that $B_{n}(X)=Z_{n}(X)=0$ and
$B_{n-1}(X)=\partial[w, C_{n-1}(Y)]$, we can choose basis for
$\tilde{b}_{n-1}(X)$ as follows:
$$\tilde{b}_{n-1}(X)=\{ [w,
b_{n-1}(Y)], [w, h_{n-1}(Y)], [w, \tilde{b}_{n-2}(Y)]\}$$ So the determinant of the
transition matrix $D_{n}$ is:
$$D_{n}=\Big[\frac{ [w, b_{n-1}(Y)], [w, h_{n-1}(Y)], [w, \tilde{b}_{n-2}(Y)]}{ [\omega, c_{n-1}(Y)]}\Big]$$
which is $D_n=A_{n-1}$ where $A_{n-1}$ denotes the corresponding determinant of the
transition matrix for $Y$.
\newline

\noindent Case 1. $n-p_n<p < n$\\
\newline
In this case, $IC_{p}(X)=C_{p}(X)$, so we only need to consider the usual chains of $X$.\\
\newline
Let $c_p(Y)=\{ \sigma_{1}^{p}(Y), \cdots, \sigma_{i_{p}}^{p}(Y)\}$
be the preferred basis of $p$-chains of $Y$.  Then \{$c_p(Y)$, $[w,
c_{p-1}(Y)]$\} is the preferred basis of $C_{p}(X)$. Choose a basis
$b_{p}(Y)=\{b_{1}^{p}(Y), \cdots, b_{k_{p}}^{p}(Y)\}$ for $B_p(Y)$,
and their lifts $\tilde{b}_{p}(Y)=\{ \tilde{b}_{1}^{p}(Y), \cdots,
\tilde{b}_{k_p}^{p}(Y)\}$, and $h_{p}(Y)=\{ h_{1}^{p}(Y), \cdots,
h_{j_{p}}^{p}(Y)\}$ the basis for $H_{p}(Y)$. Then by the fact that
$B_{p}(X)=Z_{p}(X)=\partial[w, C_{p}(X)]$, we can choose a basis for
$B_{p}(X)$ as follows:
\begin{align*}
b_{p}(X)&=\{ \partial[w, b_{p}(Y)], \partial [w, h_{p}(Y)],  \partial[w, \tilde {b}_{p-1}(Y)]  \}\\
       &=\{ b_{p}(Y),  h_{p}(Y),  \tilde{b}_{p-1}(Y)-[w, b_{p-1}(Y)]
       \}
       \end{align*}
and the lifts $\tilde{b}_{p-1}(X)$ (even though the lifts $\tilde{b}_{p-1}(X)$ depends on $B_{p-1}(X) \subset IC_{p-1}(X)$
which may not be $C_{p-1}(X)$ when $p-1=n-p_n$, $B_{p-1}(X)$ still consists of the ordinary boundaries):
$$\tilde{b}_{p-1}(X)=\{ [w,
b_{p-1}(Y)], [w, h_{p-1}(Y)], [w, \tilde{b}_{p-2}(Y)] \}$$ So the determinant of the
transition matrix $D_{p}$ is:
$$D_{p}=\Big[\frac{b_{p}(Y), h_{p}(Y), \tilde{b}_{p-1}(Y), [\omega, b_{p-1}(Y)], [\omega, h_{p-1}(Y)], [\omega, \tilde{b}_{p-2}(Y)]}{c_{p}(Y), [\omega, c_{p-1}(Y)]}\Big]$$
which is
\begin{align}\label{eqn:4.2}  D_p=A_p A_{p-1}.
         \end{align}\newline
Case 2. $p=n-p_n$.\\
\newline
In this case, we still have $IH_{p}(X)=H_{p}(X)=0$. By Lemma
\ref{icc}, $IC_{p}(X)=\{\ \eta \in C_p(X) \ | \ \partial \eta \in
C_{p-1}(Y)\ \}$.

For $\eta\in IC_{p}(X)$, write $\eta=C_{p}(Y)+[w, D_{p-1}(Y)]$. Then
$\partial \eta=\partial C_{p}(Y)+D_{p-1}(Y)-[w, \partial
D_{p-1}(Y)]$. Thus we must have $\partial D_{p-1}(Y)=0$ for
$\partial \eta \in C_{p-1}(Y)$, which implies that $D_{p-1}(Y) \in
B_{p-1}(Y)\oplus H_{p-1}(Y)$.

Thus,
\begin{align*}
\partial{(IC_{p}(X))}&=\partial (C_{p}(Y))\oplus \partial[w, B_{p-1}(Y)\oplus H_{p-1}(Y)]\\
&=\partial (C_{p}(Y))\oplus B_{p-1}(Y)\oplus H_{p-1}(Y)\\
&=B_{p-1}(Y)\oplus H_{p-1}(Y)
\end{align*}
Hence we can take $\{ [w, b_{p-1}(Y)], [w, h_{p-1}(Y)] \}$ as basis of $\tilde {B}_{p-1}(X)$. \\
\newline
The fact that $IH_{p}(X)=H_{p}(X)=0$ implies
$$IC_{p}(X)=B_{p}(X)\oplus \tilde
{B}_{p-1}(X).$$ Then the determinant of the transition matrix $D_{p}$ is
\begin{align*}
D_{p}&=\Big[\frac{b_{p}(Y), h_{p}(Y), \tilde {b}_{p-1}(Y), [\omega,
b_{p-1}(Y)],[\omega, h_{p-1}(Y)]} {c_{p}(Y), [\omega,
b_{p-1}(Y)],[\omega, h_{p-1}(Y)]}\Big]
\end{align*}
which yields \begin{align}\label{eqn:4.3} D_p=A_p
         \end{align}\newline
Case 3: $p=n-p_n-1$\\
\newline
In this case
\begin{align*}
&IC_{p}(X)=C_{p}(Y)\\
&IH_{p}(X)=\mathrm{Im} \big(H_{p}(Y)\rightarrow H_{p}(X)\big)=0
\end{align*}
Consider the following sequence:
\begin{align*}
\overset{\partial} {\cdots \longrightarrow } IC_{p+1}(X)\overset
{\partial}\longrightarrow IC_{p}(X)\overset {\partial
}\longrightarrow IC_{p-1}(X)\overset {\partial}\longrightarrow
\cdots .
\end{align*}
Then as before,
$$ \partial (IC_{p+1}(X))=B_{p}(Y)\oplus H_{p}(Y),$$
and $$ \partial [IC_{p}(X)]=\partial[C_{p}(Y)]=B_{p-1}.$$ Thus the determinant of
the transition matrix is:
\begin{align}\label{eqn:4.4}
D_{p}=\Big[\frac{b_{p}(Y), H_{p}(Y),
\tilde{b}_{p-1}(Y)}{C_{p}(Y)}\Big]=A_{p}
\end{align}
Case 4: $p<n-p_n-1$\\
\newline
In this case, $IC_{p}(X)=C_{p}(Y)$, $IH_{p}(X)=H_{p}(Y)$. Since the choice of the preferred basis on $IH_{p}(X)$ 
is compactible with that of $H_{p}(Y)$,  it is easy to see that $D_{p}=A_{p}$ \\
\newline
Combining the above results, we have:
\begin{align}\label{eqn:4.5}
\tau (IC)&=\overset {n}{\underset {p=0}\prod}(D_{p})^{(-1)^{p}}\notag\\
         &=\overset {n-p_n}{\underset {p=0}\prod}(A_{p})^{(-1)^{p}}
        \cdot \overset {n-1}{\underset {p=n-p_n+1}\prod}(A_{p} \cdot A_{p-1})^{(-1)^{p}} \cdot (A_{n-1})^{(-1)^{n}} \notag\\
        &=\overset {n-p_n-1}{\underset {p=0}\prod}(A_{p})^{(-1)^{p}}
\end{align}

Thus,
\begin{align}\label{eqn:4.6}
\ln I\tau^{\bar{p}}(X)= \ln \tau(IC)&=\overset {n-p_n-1}{\underset {p=0}\sum}(-1)^{p} \ln A_{p}\notag\\
         &=\sum_{p=0}^{n-p_n-1} (-1)^{p+1}p\,
\ln \det \Delta_p^{(c)} + (n -p_n) \sum_{p-p_n}^{n-1}
(-1)^{p+1}\ln \det \Delta_p^{(c)}.
\end{align}
Here we have used the equation $\ln A_p = -\frac{1}{2}
\sum_{k=p}^{n-1} (-1)^{k-p} \ln \det \Delta_k^{(c)}$ \cite{RS}. \qed

In Theorem \ref{mt} the basis of $IH_i^{\bar{p}}(X)$ is chosen to be compatible with those of  $H_{p}(Y)$ with respect to  (\ref{icoc}).
As discussed in the previous section, the convention for Reidemeister torsion is that the choice of the cohomology basis is determined by the Hodge theory. Thus we will now consider the issue briefly.

Recall that the conical metric on  $X=C(Y)=(0, 1]\times Y$ is $dr^2+ r^2g$, where $g$ is a metric on $Y$. There is a natural decomposition of smooth $i$-forms on $X$:
\be \label{ndof}
\alpha_i= g(r) \phi_i + f(r)   dr\wedge\omega_{i-1},
\ee
where $\phi_i,  \omega_{i-1}$ are $i$-forms resp. $(i-1)$-forms on $Y$. Following the notations from \cite{c2}, we will decorate operators on $Y$ with ``$\sim$". For example $d$ will denote the exterior derivative on $X$, while $\tilde{d}$ will denote the exterior derivative on $Y$. Similarly for the adjoints $\delta$, $\tilde{\delta}$, and the Hodge Laplacian $\Delta$, $\tilde{\Delta}$.

The following formula for the action of Laplacian is from \cite{c2}[(3.8)].
\begin{eqnarray} \label{aol}
\Delta\alpha_i & = & (-g'' - [n-2i-1]r^{-1}g') \phi_i + r^{-2}g\tilde{\Delta}  \phi_i - 2r^{-3}g dr\wedge \tilde{\delta}\phi_i   \\
&  + & dr \wedge (-f''  - [n-2i+1]r^{-1}f' + [n-2i+1]r^{-2}f )\omega_{i-1} \nonumber \\
& + &  r^{-2}f dr \wedge\tilde{\Delta} \omega_{i-1}  - 2 r^{-1}f \tilde{d} \omega_{i-1} . \nonumber
\end{eqnarray}
In particular, when $g(r)=1, f(r)=0$ and $\phi$ is harmonic on $Y$, $\alpha_i= \phi_i $ is also harmonic on $X$. Since $H^i_{(2)}(X)$ (with absolute boundary condition) is isomorphic to $(IH_i^{\bar{p}}(X))^*$ where the perversity $\bar{p}$ is taken to be the lower middle perversity, we see that for $i< n/2$ the isomorphism $H^i(Y) \simeq H^i_{(2)}(X)$ is given by $\phi_i  \rightarrow \phi_i $ for harmonic forms $\phi_i $ on $Y$.
This isomorphism is not an isometry. Indeed $\|\phi_i\|_{L^2(X)}^2=\frac{1}{n-2i} \|\phi_i\|_{L^2(Y)}^2$ ($i< n/2$). Correcting this, say via (\ref{bcfrt2}), we have (for the lower middle perversity)
\ban \ln I\tau(X) & = & \sum_{p<n/2} \frac{(-1)^p\, \ln (n-2p) }{2} b_p(Y)
 +  \sum_{p<n/2} (-1)^{p+1}p\,
\ln \det \Delta_p^{(c)}  \\
& & + (n -[\frac{n-2}{2}]) \sum_{p\geq n/2}
(-1)^{p+1}\ln \det \Delta_p^{(c)}. \ean
Here $b_p(Y)$ denotes the $p$-th Betti number of $Y$. In particular, when $n=3$ we have
$$ \ln I\tau(X) = \frac{\ln 3}{2} + \sum_{p=0}^2 (-1)^{p+1}p\,
\ln \det \Delta_p^{(c)} . $$

\section{An analytic analogue}

Following Ray-Singer's idea of defining analytic torsion as a formal
analog of the R-torsion on closed manifolds, we now study the formal
analytic analog of the intersection R-torsion (1.2), which is
intrinsic to the even dimensional cross section. That is, by
replacing the combinational Laplacian by the Hodge Laplacian, we
define an analytic invariant for an even dimensional closed
manifold.

More precisely, let $Y$ be an even dimensional closed manifold with
$m=\dim Y$. Let $p$ be an integer such that $0\leq p \leq m-1$ ($p$
corresponds to $p_n$ which is determined by a given perversity).
Given an orthogonal representation $\rho: \, \pi_1(Y)
\longrightarrow O(N)$, one has an associated flat vector bundle
$E_{\rho}$ with compatible metric on $Y$. Let $\Delta_k$ be the
Laplacian acting on differential $k$ forms on $Y$ with coefficients
in $E_{\rho}$. Then we define

\be \label{aioedm} \ln T_{p}(Y, \rho) = \frac{1}{2} \left[
\sum_{k=0}^{m-p} (-1)^{k+1} k\, \ln \det(\Delta_k) +
(m-p)\sum_{k=m-p+1}^m (-1)^{k+1} \ln \det(\Delta_k) \right] . \ee

For $p=0$, which corresponds to the minimum perversity, \[ \ln T_0(Y,
\rho) = \frac{1}{2} \sum_{k=0}^{m} (-1)^{k+1} k\, \ln
\det(\Delta_k)=0
\] is the usual analytic torsion which is trivial for even
dimensional manifolds. On the other hand, for $p=m-1$ corresponding
to the maximum perversity, \[ \ln T_{m-1}(Y, \rho) = \frac{1}{2}
\sum_{k=1}^{m} (-1)^{k+1} \ln \det(\Delta_k).
\] The more interesting cases are given by $p=\frac{m}{2}-1$ and
$p=\frac{m}{2}$ corresponding to the lower and upper middle
perversity, respectively. In these cases, we have
\ban \ln T_{\frac{m}{2}-1}(Y, \rho) & = & \frac{1}{2} \left[
\sum_{k=0}^{\frac{m}{2}+1} (-1)^{k+1} k \ln \det(\Delta_k) +
(\frac{m}{2} + 1)\sum_{k=\frac{m}{2}+2}^{m} (-1)^{k+1} \ln
\det(\Delta_k) \right] \\
& = & \frac{1}{2} \left[ \sum_{k=0}^{\frac{m}{2}} (-1)^{k+1} k \ln
\det(\Delta_k) + (\frac{m}{2} + 1)\sum_{k=\frac{m}{2}+1}^{m}
(-1)^{k+1} \ln \det(\Delta_k) \right] \ean
and
\[ \ln T_{\frac{m}{2}}(Y, \rho) = \frac{1}{2} \left[
\sum_{k=0}^{\frac{m}{2}} (-1)^{k+1} k \ln \det(\Delta_k) +
\frac{m}{2}\sum_{k=\frac{m}{2}+1}^{m} (-1)^{k+1} \ln \det(\Delta_k)
\right] . \]
 When $Y$ is oriented, we can actually use
Poincare duality to write it in terms of the Laplacians on half of
the degrees. For example, for $p=\frac{m}{2}-1$ corresponding to the
lower middle perversity, we have
 \be \label{aioedm2} \ln T_{\frac{m}{2}-1}(Y, \rho) = \frac{1}{2} \left[
\sum_{k=0}^{\frac{m}{2}-1} (-1)^{k+1} (k+ \frac{m}{2} + 1) \, \ln
\det(\Delta_k) + (-1)^{\frac{m}{2}+1} \frac{m}{2} \ln
\det(\Delta_{\frac{m}{2}}) \right] . \ee
These bear close resemblance to the so called Cheeger's half torsion \cite{}.

To investigate what kind of invariant $\ln T_{p}(Y, \rho)$ defines,
we now look at its variation under metric change. Let $g(u)$ be a
family of Riemannian metrics on $Y$ and $\Delta_k(u)$ the
corresponding Laplacian (when there is no ambiguity we will often
write $\Delta_k$ instead of $\Delta_k(u)$.  Let $\dot{\star}=
d\star/du$ and $\alpha=\star^{-1} \dot{\star}$.  Denote by
$E_k(t)=e^{-t\Delta_k(u)}$ the heat kernel and let $E_k=E_k^{ex} +
E_k^{ce} + E_k^h$ denote the Hodge decomposition of $E_k$ into its
exact, coexact and harmonic parts. We have the following result
regarding the variation of $\ln T_{p}(Y, \rho)$.

\begin{theo} \label{vfai} The variation of $\ln T_p(Y, \rho)$ is given by
\ban \frac{d}{du} \ln T_p(Y, \rho) & = & \hf \sum_{k=0}^{m-p-1}
(-1)^{k+1} {\rm Tr}(P_{H^k} \alpha) + \hf \sum_{k=0}^{m-p-1}
(-1)^{k+1} {\rm LIM}_{t\rightarrow 0} {\rm
Tr}(e^{-t\Delta_k} \alpha)\\
& + & (-1)^{m-p+1} \hf {\rm LIM}_{t\rightarrow 0} {\rm
Tr}(E_{m-p}^{ex}(t) \alpha), \ean where $P_{H^k}$ denote the
projection onto the cohomology $H^k$ and ${\rm LIM}_{t\rightarrow 0}
{\rm Tr}(E_{m-p}^{ex}(t) \alpha)$ denotes the constant term in the
asymptotic expansion of ${\rm Tr}(E_{m-p}^{ex}(t) \alpha)$.
\end{theo}

Before we give the proof of our theorem, we need the following
result from \cite{c1} (compare also with \cite{RS}) concerning the
variation of heat kernel.

\begin{theo}[Cheeger] The variation of the trace of the heat kernel $E_k$ is given by
\ban \frac{d}{du} \tr (E_k(t)) & = & -t[ \tr
\left(\Delta_{k+1}E_{k+1}^{ex} \alpha \right) - \tr
\left(\Delta_{k}E_{k}^{ce} \alpha \right)
+  \tr \left(\Delta_{k}E_{k}^{ex} \alpha \right) - \tr \left(\Delta_{k-1}E_{k-1}^{ce} \alpha \right) ]\\
& = & t \frac{d}{dt} [ \tr \left(E_{k+1}^{ex} \alpha \right) - \tr
\left(E_{k}^{ce} \alpha \right) + \tr \left(E_{k}^{ex} \alpha
\right) - \tr \left(E_{k-1}^{ce} \alpha \right) ]. \ean
\end{theo}

The following lemma is an immediate consequence of Cheeger's result.

\begin{lem} For any integer $q$, $0\leq q \leq m$, we have
\be \label{varf1} \frac{\pt}{\pt u} \sum_{k=0}^q (-1)^k k \, {\rm tr} (E_k(t))= t
\frac{\pt}{\pt t} [ \sum_{k=0}^q (-1)^k {\rm tr} (\underline{E}_k(t) \alpha) +
(-1)^q q {\rm tr} (E_{q+1}^{ex}(t)\alpha) +  (-1)^{q+1} (q+1) {\rm tr}
(E_{q}^{ce}(t)\alpha)]. \ee
Similarly, for any integer $r$, $0 \leq r \leq
m$, \be \label{varf2} \frac{\pt}{\pt u} \sum_{k=r}^m (-1)^k  {\rm tr} (E_k(t))= t
\frac{\pt}{\pt t} [ (-1)^r {\rm tr} (E_{r}^{ex}(t)\alpha) + (-1)^{r-1}
{\rm tr} (E_{r-1}^{ce}(t)\alpha)]. \ee

\end{lem}

With these results at our disposal, we are now ready to prove the variational formula for our analytic invariant.
\newline

{\em Proof of Theorem \ref{vfai}}: Define for $\Re \, s$
sufficiently large
\[ f(u, s)= \hf \big[ \sum_{k=0}^{m-p} (-1)^k k \int_0^{\infty} t^{s-1} {\rm Tr}(e^{-t[\Delta_k + P_{H^k}]})\, dt
+ (m-p) \sum_{k=m-p+1}^m (-1)^k \int_0^{\infty} t^{s-1} {\rm
Tr}(e^{-t[\Delta_k + P_{H^k}]})\, dt \big].\] Then $f(u, s)$ has a
meromorphic extension to the whole complex $s$-plane with a simple
pole at $s=0$. Indeed, since
\[ {\rm Tr}(e^{-t[\Delta_k + P_{H^k}]})= {\rm
Tr}(e^{-t\underline{\Delta}_k}) + e^{-t} \dim H^k, \] we have
\[ {\rm Res}_{s=0} f(u, s)= \hf \big[ \sum_{k=0}^{m-p} (-1)^k k
A_{m/2, k} + (m-p) \sum_{k=m-p+1}^m (-1)^k A_{m/2, k} \big], \]
where $A_{m/2, k}$ denotes the constant term in the asymptotic
expansion of ${\rm Tr}(e^{-t \Delta_k})$. Now let
\[ \tilde{f}(u, s)= f(u, s) - \Gamma(s){\rm Res}_{s=0} f(u, s). \]
Then $\tilde{f}$ is holomorphic at $s=0$ and we have
\[ \tilde{f}(u, 0)=\ln T_p(Y, \rho). \]

Now, for $\Re \, s$ sufficiently large
\ban \frac{\pt}{\pt u}\int_0^{\infty} t^{s-1} {\rm Tr}(e^{-t[\Delta_k + P_{H^k}]}) \, dt & =  & \int_0^{\infty} t^{s-1}  \frac{\pt}{\pt u}
{\rm Tr}(e^{-t[\Delta_k + P_{H^k}]}) \, dt  \\
& =  &  \int_0^{\infty} t^{s-1}  \frac{\pt}{\pt u}
{\rm Tr}(e^{-t\Delta_k }) \, dt  \\
\ean
Hence, using (\ref{varf1}), (\ref{varf2}), we derive
\ban  \frac{\pt}{\pt u} f(u, s) & = & \hf \big[ \sum_{k=0}^{m-p} (-1)^k  \int_0^{\infty} t^{s} \frac{\pt}{\pt t}{\rm Tr}(\underline{E}_k(t)\alpha)\, dt
+  (-1)^{m-p+1} \int_0^{\infty} t^{s} \frac{\pt}{\pt t} {\rm
Tr}(E^{ce}_{m-p}(t)\alpha)\, dt \big] \\
& = & s  \hf \big[ \sum_{k=0}^{m-p} (-1)^{k+1}  \int_0^{\infty} t^{s-1}{\rm Tr}(\underline{E}_k(t)\alpha)\, dt
+  (-1)^{m-p} \int_0^{\infty} t^{s-1}  {\rm
Tr}(E^{ce}_{m-p}(t)\alpha)\, dt \big] \\
& = & s  \hf \big[ \sum_{k=0}^{m-p-1} (-1)^{k+1}  \int_0^{\infty} t^{s-1}{\rm Tr}(\underline{E}_k(t)\alpha)\, dt
+  (-1)^{m-p+1} \int_0^{\infty} t^{s-1}  {\rm
Tr}(E^{ex}_{m-p}(t)\alpha)\, dt \big]
\ean

It follows then that \ban \frac{\pt}{\pt u}\ln T_p(Y, \rho) & = &
\hf \sum_{k=0}^{m-p} (-1)^k {\rm Tr}(P_{H^k} \alpha) + \hf
\sum_{k=0}^{m-p-1} (-1)^{k+1} {\rm LIM}_{t\rightarrow 0} {\rm
Tr}(e^{-t\Delta_k} \alpha)\\
& + & (-1)^{m-p+1}\hf {\rm LIM}_{t\rightarrow 0} {\rm
Tr}(E_{m-p}^{ex}(t) \alpha). \ean \qed

Just like Cheeger's half torsion, we note the special property of the invariant for $p=m/2$ under conformal change.

\begin{coro} Under a family of conformal  changes, for $p=m/2$, 
the variation of $\ln T_{\frac{m}{2}}(Y, \rho)$ is local in the sense that
\[ \frac{d}{du} \ln T_{\frac{m}{2}}(Y, \rho)  =  \hf \sum_{k=0}^{\frac{m}{2}-1}
(-1)^{k+1} {\rm Tr}(P_{H^k} \alpha) + \hf \sum_{k=0}^{\frac{m}{2}-1}
(-1)^{k+1} {\rm LIM}_{t\rightarrow 0} {\rm
Tr}(e^{-t\Delta_k} \alpha). \]
\end{coro}

{\em Proof.} If $g=e^{2f}g_0$ is a conformal change of $g_0$, then its Hodge star on the $p$-forms is given by
$*_g=e^{(2p-m)f} *_0$ in terms of the star operator of $g_0$. It follows that $\alpha=0$ on $m/2$-forms. Our result follows 
from Theorem \ref{vfai}. \qed

\section{R-torsion of the Mayer-Vietoris sequences}

Consider an $(m+1)$-dimensional Riemannian manifold $X$ with
isolated conical singularity. Thus, $X=C(Y)\cup M$, where $M$ is a
compact manifold with boundary and $\partial M=Y$. It is understood
in this section that the collar neighborhoods of the boundaries of
$M$ and $C(Y)$ are extended so that they form an open cover of $X$.
We assume that $m+1$ is odd.

As we mentioned, the general Mayer-Vietoris Principle reduces the
torsion of $X$ to that of $C(Y)$, $M$ as well as the torsion of the
Mayer-Vietoris sequence in the intersection cohomology. We now
examine the torsion of the Mayer-Vietoris sequence.

We use the $L^2$-cohomology interpretation of the intersection
cohomology in this setting \cite{c2}. The Mayer-Vietoris sequence
goes
\begin{align} \label{mvs}
\cdots \longrightarrow H_{(2)}^{q}(Y) \overset {d^*}\longrightarrow
H_{(2)}^{q+1}(X)\longrightarrow H_{(2)}^{q+1}(M)\oplus
H_{(2)}^{q+1}(C(Y))\longrightarrow H_{(2)}^{q+1}(Y) \longrightarrow
\cdots .
\end{align}

First, we have the following

\begin{lem} For the Mayer-Vietoris long exact
sequence in cohomology (\ref{mvs}), \\
a). its part for $q\leq m/2$ splits
into the following short exact sequences:
\begin{align}\label{ses}
0\longrightarrow H_{(2)}^{q}(X)\longrightarrow H_{(2)}^{q}(M)\oplus
H_{(2)}^{q}(C(Y))\longrightarrow H_{(2)}^q(Y)\longrightarrow 0
\end{align}
b). further,
\begin{align}\label{eqn:1.5} 0\longrightarrow
H_{(2)}^{q}(X)\longrightarrow H_{(2)}^{q}(M)\oplus
H_{(2)}^{q}(C(Y))\longrightarrow H_{(2)}(Y)\longrightarrow 0
\end{align} is a split short exact sequence. \\
c). the part of the Mayer-Vietoris sequence for $q> m/2$ is
naturally isomorphic to the truncated exact sequence for the pair
$(M, Y)$:
\begin{align}
 H^{m/2}(Y)\longrightarrow H^{m/2+1}(M,Y)\longrightarrow H^{m/2+1}(M)\longrightarrow H^{m/2+1}(Y) \longrightarrow \cdots \longrightarrow H^m(Y) .\notag \\
\end{align}
 \end{lem}

\Pf  For a). we only need to show that, when $q\leq m/2$, $Im
 (d^*)=0$. Let $\rho_1, \rho_2$ be a partition of unity subordinate
 to the open cover of $X$ by $M, C(Y)$. That is, $\rho_1, \rho_2 \in
 C^{\infty}(X)$, $0 \leq \rho_1, \rho_2 \leq 1$, $\rho_1+ \rho_2=1$
 and ${\rm supp} \, \rho_1 \subset M, \ {\rm supp} \, \rho_2 \subset
 C(Y)$. Then, for a closed $q$-form on $Y$,
 \begin{align*}
d^*[w]=\begin{cases}
[-d(\rho_2 w)]  & on  \quad M ,\\
[d(\rho_1w)] & on \quad C(Y).
\end{cases}
\end{align*}
Here $w$ is extended trivially along radial directions hence defines
a $q$-form in a collared neighborhood of $Y$ in $X$. In fact,
$d^*[w]$ is supported in this collared neighborhood and, interpreted
properly, either $[-d(\rho_2 w)]$ or $[d(\rho_1w)]$ defines
$d^*[w]$. Now, $d^*[w]= [-d(\rho_2 w)]$. By the result of \cite{c2},
for $q\leq m/2$, $w$ defines an $L^2$ form on $C(Y)$. This shows
that $d^*[w]$ is exact in $L^2$ cohomology. Hence $d^*[w]=0$.

The statement b). is clear since these are short exact sequences of
vector spaces. They can also be seen directly as follows. We show
that the composition $p \, i^*$ in the follwoing diagram
\begin{align*}
H_{(2)}^{q}(X)\overset{i^*}\longrightarrow H^{q}(M)&\oplus H_{(2)}^{q}(C(Y))\longrightarrow H^q(Y)\\
 &\downarrow p\\
 H&^{q}(M)
\end{align*}
is an isomorphism. Here $p$ is the projection onto the first factor.
Indeed, for any $w\in H_{(2)}^q(X)$, $p \, i^*{w}=p\, (i_{M}^{*}{w},
i_{C(Y)}^*{w})= i_{M}^{*}{w}$.  If $i_{M}^{*}{w}$ is an exact form,
$i_{M}^{*}{w}=d\eta_{2}$ then
$i_{Y}^{*}{i_{M}^*{w}}=i_{Y}^*(d\eta_2)=d(i_Y^*{\eta_2})$ is exact
on $Y$. By \cite{c2}, for $q \leq m/2$, $i_Y^*{\eta_2}$ defines an
$L^2$ form on $C(Y)$. Since the cohomology class of a closed form on
$C(Y)$ is uniquely determined by its restriction on $Y$ \cite{c2},
we see that $i_{C(Y)}^*(w)$ is exact. It follows then that
$i^{*}(w)=(i_{M}^{*}{w}, i_{C^*_{0,1}(N)}^{*}{w})$ is exact. Namely
$[i^{*}w]=0$ on $H_{2}^{*}(M)\oplus H_{2}^{*}(C^*_{0,1}(N))$. So
$[w]=0$ on $H_{(2)}^{q}(X)$ by the injectivity of the short exact
sequence. This shows that $p \, i^*$ is injective.

For the surjectivity, take $\eta \in H_{(2)}^{q}(M)$. Let
$\xi=i_{Y}^{*}(\eta)\in H_{(2)}^{q}(Y)$. Then $\xi$ extends to an
$L^2$ form on $C(Y)$ which is cohomologous with the restriction of
$\eta$ in a collared neighborhood of $Y$. It follows that $(\eta,
\xi)$ is the image of some element of $H_{(2)}^{q}(X)$, say $ w$.
then $p\, i^{*}(w)=\eta
 $

Part c). follows from the natural isomorphisms $H_{(2)}^q(X) \cong
H^q(M, Y)$, $H_{(2)}(C(Y)) \cong 0$ for $q>m/2$ \cite{c2}.

 \qed

 \begin{lem}
 For a split short exact sequence
 \[ 0 \longrightarrow V_1 \overset{i} \longrightarrow V_2 \overset{p} \longrightarrow V_3 \longrightarrow 0 \]
 with preferred bases $c_1, c_2, c_3$, its R-torsion is determined by $i(c_1), j(c_3)$ and $c_2$, where $j$ is an homomorphism
 from $V_3$ to $V_2$ such that $p j = {\rm id}$. In fact, the
 R-torsion is given by
 \[ | [ i(c_1)j(c_3)/c_2] | \]
 \end{lem}

 \Pf We choose $b_1=0$, $b_2=i(c_1)$, and $b_3=c_3$ and set
 $\tilde{b}_1=c_1$, $\tilde{b}_2=j(c_3)$ and $\tilde{b}_3=0$.
 The lemma follows.
 \qed

A split short exact sequence can be written as
 \[ 0 \longrightarrow V_1 \overset{i} \longrightarrow V_1\oplus V_3 \overset{p} \longrightarrow V_3 \longrightarrow 0, \]
 where $i$ is not necessarily the natural inclusion, nor $p$ the
 natural projection.

 \begin{lem}
 For a split short exact sequence
 \[ 0 \longrightarrow V_1 \overset{i} \longrightarrow V_1\oplus V_3 \overset{p} \longrightarrow V_3 \longrightarrow 0 \]
 with preferred bases $c_1, c_1\oplus c_3, c_3$, consider the natural projection $p_1: \, V_1\oplus V_3 \longrightarrow V_1$ onto the first factor and the
 natural inclusion $i_2:\, V_3 \longrightarrow V_1 \oplus V_3$ of the second factor.
 If $ p_1i:\, V_1 \longrightarrow V_1$ is an isometry with respect to the inner product induced by the
 preferred basis $c_1$ and $pi_2={\rm id}:\, V_3 \longrightarrow V_3$,  then the R-torsion
 of the short exact sequence is
 trivial.
 \end{lem}

\Pf Using the lemma above we just need to compare the basis
$i(c_1)\oplus c_3$ with $c_1 \oplus c_3$. Since $pi_1$ is an
isometry, we might as well replace $c_1 \oplus c_3$ with $pi_2(c_1)
\oplus c_3$. Then clearly, the transition matrix from $i(c_1)\oplus
c_3$ to $pi_2(c_1) \oplus c_3$ is an upper triangular matrix with
all diagonal entries one. The lemma follows.\qed

Combining the above results, we obtain the main result of this
section on the R-torsion of the Mayer-Vietoris sequence.
\begin{theo}
Assume that the Witt condition $H^{\frac{m}{2}}(Y)=0$ holds. Then
the R-torsion of the Mayer-Vietoris sequence in intersection
cohomology
$$ \cdots \longrightarrow IH^{q}_{(2)}(Y)\longrightarrow IH_{(2)}^{q+1}(X)\longrightarrow IH_{(2)}^{q+1}(M)\oplus
IH_{(2)}^{q+1}(C(Y))\longrightarrow IH_{(2)}^{q+1}(Y)\longrightarrow
\cdots
$$
is equal to the R-torsion of the truncated exact sequence of the
pair $(M,Y)$
$$ 0 \longrightarrow H^{\frac{m}{2}+1}(M,Y) \longrightarrow H^{\frac{m}{2}+1}(M)\longrightarrow
H^{\frac{m}{2}+1}(Y)\longrightarrow H^{\frac{m}{2}+2}(M,Y)
\longrightarrow \cdots $$
\end{theo}


\begin{thebibliography}{99}


\bibitem{bz}
J.-M. Bismut and W.~Zhang.
\newblock An extension of a theorem by {Cheeger and M\"uller}.
\newblock {\em Ast\'erisque}, 205, 1992.

\bibitem{BM} J. Br\"uning, X. Ma.
\newblock  An anomaly formula for Ray-Singer metrics on manifolds with
boundary.
\newblock {\em Geom. Funct. Anal.} 16 (2006), no. 4,
767--837.

\bibitem {BFKM} D. Burghelea, L. Friedlander, T. Kappeler, and P.
McDonald.
\newblock Analytic and Reidemeister Torsion for
Representations in Finite Type Hilbert Modules.
\newblock {\em Geom.
Funct. Anal.}, 6 (1996), pp. 751-859.

\bibitem{c1}  J. Cheeger.
\newblock  Analytic torsion and the heat equation.
\newblock {\em Ann. Math.} 109 (1979) 259-322.

\bibitem{c2}  J. Cheeger.
\newblock Spectral geometry of singular Riemannian spaces.
\newblock {\em J. Diff. Geom.}  18, 575-657 (1983).

\bibitem{c3}  J. Cheeger.
\newblock $\eta$-invariants, the adiabatic approximation and conical singularities. I. The adiabatic approximation. \newblock {\em  J. Differential Geom.}  26  (1987),  no. 1, 175--221.


\bibitem{d1}
X.~Dai.
\newblock Adiabatic limits, nonmultiplicativity of signature, and {Leray}
  spectral sequence.
\newblock {\em J. Amer. Math. Soc.}, 4:265--321, 1991.

\bibitem{d2}
X.~Dai.
\newblock APS boundary conditions, eta invariants and adiabatic limits
\newblock {\em Trans. AMS}, 354, pp. 107-122

\bibitem{DF} X. Dai, H. Fang.
\newblock  Analytic torsion and R-torsion for manifolds with boundary.
\newblock {\em Asian J. Math.}  4  (2000),  no. 3, 695--714.

\bibitem{df1}  X. Dai, D. Freed.
\newblock  Eta invariants and determinant lines.
\newblock {\em  J. Math. Phys.},  35(1994), 5155-5194.

\bibitem{df2}  X. Dai, D. Freed.
\newblock   Invariants eta et droites determinants.
\newblock {\em C. R. Acad. Sci., Paris, Series I}, 320(1995),  585-591.

\bibitem{dar}  A. Dar.
\newblock Intersection $R$-torsion and analytic torsion for pseudomanifolds.
\newblock {\em Math. Z.}  194  (1987),  no. 2, 193--216.

\bibitem {F} W. Franz.
\newblock  \"{U}ber die Torsion einer \"{u}berdeckrung.
\newblock {\em J. Reine Angew. Math.}, 173 (1935), pp. 245-254

\bibitem{gm1} M. Goresky, R. MacPherson.
\newblock  Intersection homology theory.
\newblock {\em Topology}  19  (1980), no. 2, 135--162.

\bibitem{gm2} M. Goresky, R. MacPherson.
\newblock Intersection homology. II.
\newblock {\em Invent. Math.}  72  (1983), no. 1, 77--129.

\bibitem {hs} L. Hartmann and M. Spreafico.
\newblock The analytic torsion of a cone over an odd dimensional manifold.
\newblock  {\em  J. Geom. Phys.}, 61 (2011), no. 3, 624–657.


\bibitem {K} Frances Kirwan.
\newblock {\em  An Introduction to Intersection Homology Theory}.
\newblock Pitman Reasearch Notes in Mathematics Series., 187.


\bibitem {LR} J. Lott and M. Rothenberg.
\newblock Analytic Torsion for Group Actions.
\newblock {\em J. Diff. Geom.}, 34 (1991), pp. 431-481.

\bibitem {lu}  W. L\"{u}ck.
\newblock  Analytic and Topological Torsion for Manifolds with Boundary and
Symmetry.
\newblock {\em J. Diff. Geom.}, 37 (1993), pp. 263-322.

\bibitem{M} J. Milnor.
\newblock Whitehead torsion.
\newblock {\em Bull. Amer. Math. Soc.} 72 1966 358--426.


\bibitem {mu1} W. M\"{u}ller.
\newblock  Analytic Torsion and $R$-torison of Riemannian Manifolds.
\newblock {\em Adv. in Math.}, 28 (1978), pp. 233-305.

\bibitem {mu2} W. M\"{u}ller.
\newblock Analytic Torsion and $R$-torison for Unimoduler
Representations.
\newblock {\em J. Amer. Math. Soc.}, 6 (1993), pp.
31-34.

\bibitem {mv} W. M\"{u}ller and B. Vertman
\newblock The metric anomaly of analytic torsion on manifolds with conical singularities.
\newblock {\em Comm. Partial Differential Equations}, 39 (2014), no. 1, 146–191.

\bibitem {RS} D. B. Ray and I. M. Singer.
\newblock $R$-Torsion and the Laplacian on Riemannian Manifolds.
\newblock {\em Advances in Mathematics} 7., 145-210 (1971).

\bibitem {R} K. Rediemeister.
\newblock  Homotopieringe und Linsenra$\ddot{u}$m.
\newblock in {\em Hamburger Abhandl.}, 11, 1935, pp. 102-109

\bibitem {Ro} S. Rosenberg.
\newblock  Nonlocal invariants in index theory. 
\newblock  {\em Bull. Amer. Math. Soc.}, (N.S.) 34 (1997), no. 4, 423–433. 

\bibitem {s} M.  Spreafico.
\newblock  Zeta function and regularized determinant on a disc and on a cone.
\newblock {\em  J. Geom. Phys.}, 54 (2005), no. 3, 355–371.

\bibitem {T} Vladimir Turaev.
\newblock {\em Introduction to Combinational Torsions}.
\newblock Lectures in Mathematics., 2000.

\bibitem {ve} B. Vertman.
\newblock  Analytic torsion of a bounded generalized cone.
\newblock  {\em Comm. Math. Phys.}, 290 (2009), no. 3, 813–860.

\bibitem {V} S. M. Vishik.
\newblock  Generalized Ray-Singer Conjecture I: a Manifold with a Smooth
Boundary.
\newblock {\em Comm. Math. Phys.}, 167 (1995), pp. 1-102.



\end{thebibliography}
\end{document}